\documentclass[]{interact}

\usepackage[caption=false]{subfig}
\usepackage[natbibapa,nodoi]{apacite}
\usepackage{amsfonts,bm,amssymb}
\usepackage{array}

\usepackage{graphics} 
 \usepackage{epsfig}
 \usepackage{hhline}
 \usepackage{resizegather}
 \usepackage{graphicx, epsfig,color}
 \usepackage{enumerate}
\usepackage{mathtools}
\usepackage{tabularx}
\usepackage{makecell}
\usepackage{epsfig}
\usepackage{algorithm2e}
\usepackage{algorithmic}
\usepackage{resizegather}
\usepackage{url}

\usepackage{cases}
\usepackage{physics}
\usepackage{diagbox,multirow} 
\usepackage{siunitx}
\usepackage{siunitx}

\usepackage{graphicx} 
\usepackage{booktabs} 
\usepackage{xparse}   
\NewDocumentCommand{\rot}{O{45} O{1em} m}{\makebox[#2][l]{\rotatebox{#1}{#3}}}%

\usepackage{diagbox}
\usepackage{rotating}
\usepackage{graphicx}
\usepackage{makecell, caption, booktabs}
\DeclareCaptionLabelSeparator{space}{ } 

\usepackage{diagbox}
\usepackage{lscape}
\newcolumntype{s}{>{\hsize=.83\hsize}X}

\newcolumntype{L}{>{\hsize=1.4\hsize}X}
\newcolumntype{M}{>{\hsize=.3\hsize}X}
\newcolumntype{E}{>{\hsize=1\hsize}X}
 
\newcolumntype{D}{>{\hsize=1.7\hsize}X}
\newcolumntype{F}{>{\hsize=0.8\hsize}X}

\DeclarePairedDelimiter\floor{\lfloor}{\rfloor}

\newcommand{\RNum}[1]{\uppercase\expandafter{\romannumeral #1\relax}}

\theoremstyle{remark}
\newtheorem{remark}{Remark}
\theoremstyle{plain}

\newtheorem{assumption}{Assumption}

\ifdefined\pdfminorversion
\fi

\begin{document}

\renewcommand{\articletypefont}{\normalfont\fontsize{12}{14}\selectfont}
\articletype{{\fontsize{9}{11}\selectfont\itshape This is the accepted version of: Dabiri, A., \"{O}n\"{u}r, G., Gros, S., \& De Schutter, B. (2026). Optimization-based coordination of traffic lights and automated vehicles at intersections. Journal of Intelligent Transportation Systems: Technology, Planning, and Operations.}\\[9pt]RESEARCH ARTICLE}

\title{Optimization-based Coordination of Traffic Lights and Automated Vehicles at Intersections}

\author{
\name{Azita Dabiri\textsuperscript{a}, Giray \"{O}n\"{u}r\textsuperscript{a}\thanks{CONTACT Giray \"{O}n\"{u}r. Email: g.oenuer@tudelft.nl}, Sebastien Gros\textsuperscript{b} and Bart De Schutter\textsuperscript{a}}
\affil{\textsuperscript{a}Delft Center for Systems and Control, Delft University of Technology, Mekelweg 2, 2628 CD Delft, The Netherlands; \textsuperscript{b}Department of Engineering Cybernetics, Norwegian University of Science and Technology, NO-7491 Trondheim, Norway}
}

\maketitle

\begin{abstract}
This paper tackles the challenge of coordinating traffic lights and automated vehicles at signalized intersections, formulated as a constrained finite-horizon optimal control problem. The problem falls into the category of mixed-integer nonlinear programming, posing challenges for solving large instances. To address this, we introduce a decomposition approach consisting of an upper-level problem for traffic light timing allocation and a set of lower-level problems that generate appropriate commands for automated vehicles in each intersection movement. By leveraging solutions from the lower-level problems and employing parametric optimization techniques, we solve the upper-level problem using a standard sequential quadratic programming approach. The paper concludes by presenting an illustrative numerical example that highlights the effectiveness of our algorithm compared to scenarios where no coordination between traffic lights and vehicles exists.
\end{abstract}

\begin{keywords}
Receding Horizon Control; Parametric Optimization; Connected and Autonomous Vehicles; Intelligent Transportation Systems
\end{keywords}

\section{Introduction}
In recent years, there have been significant advancements in communication systems, sensing technology, driver assistance systems, and automated driving technologies. These developments, coupled with methodological improvements in fields like control theory and numerical optimization, have opened up new opportunities to address the challenges of urban traffic management. Consequently, this domain has become one of the most attractive application areas for research in cyber-physical systems.

Before delving into the contributions of this work, in the remainder of this section, we provide an overview of existing approaches in the literature that aim to mitigate congestion at intersections, enhance vehicle energy efficiency, and reduce emissions. While this section is not intended to offer an exhaustive literature survey, it aims to present a snapshot of research directions that all share a common objective: achieving efficient urban mobility. These approaches are categorized into four main areas: 1) methods for efficient traffic light control at intersections, 2) methods for optimizing the trajectories of vehicles approaching signalized intersections, 3) methods for coordinating vehicles approaching intersections without traffic lights, and 4) methods for the joint optimization of traffic signal timing and vehicle trajectories. For a comprehensive literature survey, interested readers are referred to \citet{GUOsurvey}.
\subsection{Traffic light control at intersections}
Control of traffic lights at intersections is a conventional approach for urban traffic control. The primary goal is to develop control strategies for signal timing with performance objectives such as reducing travel time, emissions, total queue length, or a combination of these factors. Approaches in the literature vary based on several factors, including the control strategy's nature, whether it is local~\citep{De-Schutter:1998aa, Haddad:2010aa}, centralized~\citep{Weg:2019aa}, or distributed~\citep{Islam:2017aa, Zhou:2015aa}, as well as whether they are model-based~\citep{Jamshidnejad:2018aa, Kulcsar:2015aa, Weg:2019aa, Islam:2017aa, Zhou:2015aa} or model-free~\citep{El-Tantawy:2013aa, Bazzan:2010aa, Zaidi:2016aa}.

Additionally, these approaches may employ various control methods, such as fuzzy logic~\citep{Rahman:2009aa}, model predictive control~\citep{Jamshidnejad:2018aa, Zhou:2015aa, Weg:2019aa}, and reinforcement learning~\citep{El-Tantawy:2013aa, Bazzan:2010aa}. They also utilize different traffic models, including the S-model~\citep{Jamshidnejad:2018aa, Lin:2011aa}, the link transmission model~\citep{Weg:2019aa, Hajiahmadi:2016aa}, and the cell transmission model~\citep{Lo:1999aa, Mohebifard:2019aa}.
It is important to note that many of these approaches incorporate aggregated variables such as traffic flow, traffic density, and average speed into their formulations.
\subsection{Vehicular trajectory control}
Acceleration and deceleration resulting from drivers' reactions to traffic lights, as well as stop-and-go movements, lead to increased fuel consumption and emissions. With the assumption of having exact and accurate information on the phase and timing of the traffic lights, the literature offers various approaches to control vehicle trajectories effectively and minimize these economic and environmental costs. These approaches differ in control methods, control objectives, types of vehicles considered, the level of vehicle connectivity, and vehicle autonomy. For single-vehicle scenarios, methods such as model predictive control~\citep{ecodrive1, ecodrive-fastMPC}, dynamic programming~\citep{ecodrive3, ecodrive4}, and other optimization-based algorithms~\citep{ecodrive2, ecodrive-queue} have been used. With the growing market penetration of electric and hybrid electric vehicles, studies focusing on their energy efficiency in signalized urban areas have gained significant attention as well~\citep{Luo2017, Yang2019, Yu2015}.

Recent studies have further explored model predictive control (MPC) and Pontryagin's minimum principle (PMP) for improving vehicle energy efficiency. MPC addresses the control problem over a moving prediction horizon, while PMP provides necessary conditions for optimality that can be used to solve the associated optimal control problem. For instance, in \citet{Guo2024PMP}, PMP is employed to solve a predictive energy management problem for fuel cell plug-in hybrid electric vehicles, using costate bounds to facilitate the online solution. In \citet{Xu2024EcoDriving}, an MPC-based eco-driving strategy is developed for battery electric vehicles with a declutchable powertrain, exploiting free-coasting to reduce energy consumption.

Although control of a vehicle trajectory has been shown to reduce energy consumption or travel time effectively, vehicular trajectory control is expected to be more effective when applied to a coordinated vehicle group. Research has explored platoon-based strategies for cooperation between human-driven and automated vehicles. For instance, in \citet{Zhao2018}, model predictive control is employed to minimize fuel consumption within platoons. In \citet{Liu2016}, platoons are reorganized to maximize the number of vehicles passing the intersection during green phases, considering safety, passenger comfort, and fuel consumption. Although promising results in reducing energy consumption or travel time are obtained, these approaches rely on the availability of accurate traffic light phase and timing information without control over how the traffic signal evolves.
\subsection{Coordination of connected vehicles at intersections without traffic lights}
Advancements in automated driving and connected vehicles have not only enabled researchers to envision innovative solutions~\citep{Rios-Torres:2017aa} but have also facilitated field-test validation~\citep{Hult:2019aa}. In these tests, coordination between vehicles plays the role of traffic lights, determining the right of way while avoiding collisions at intersections. This coordination problem can be framed as a constrained optimal control problem, involving finding the (sub)optimal order for vehicles to cross and their optimal trajectories, often with the aim of minimizing metrics like fuel consumption or travel time. However, due to its complexity, the problem is simplified in \citet{Zhang:2016aa} and \citet{Riegger:2016aa} by fixing the crossing order using heuristics such as a first-come-first-served policy and finding the optimal vehicle trajectories given this fixed order. A few contributions have been proposed that optimize both the trajectories and the crossing order of vehicles by formulating the problem in a mixed-integer programming scheme~\citep{Bali:2018aa, Hult:2018aa}. Despite their promising outcomes, these approaches do not adequately address the safety needs of vulnerable road users, such as pedestrians and cyclists, at intersections without traffic lights.
\subsection{Vehicle trajectories and traffic light control  }\label{sec:vehiclelight}
Research into the joint control of vehicle trajectories and signal timing has gained traction only in recent years, resulting in a limited number of studies within this category. One of the early contributions to this field can be found in \citet{LI:2014}, where the joint optimization of signal timing and automated vehicle trajectories was explored. To simplify the problem, the study introduced a four-component trajectory with constant acceleration in each part and used enumeration techniques to find optimal signal timing and trajectories in a straightforward intersection scenario featuring two single-lane through approaches.
This enumeration-based approach, while insightful, is not applicable to more complex intersections and is not considered optimal due to the simplifications made in trajectory generation.

In \citet{Biao_xu-2019}, a cooperative control method tackles the complexity of a 4-leg intersection using a bi-level framework. In the upper level, the optimization focuses on traffic signal timing and vehicle arrival times to minimize the total travel time of all vehicles. Meanwhile, in the lower level, vehicles optimize their trajectories to minimize fuel consumption while adhering to the scheduled times from the upper level.
However, this framework assumes a fixed cycle length and relies on enumerating all possible signal timings, which becomes computationally demanding with longer cycles or smaller time discretization intervals. Additionally, because the upper-level optimization lacks explicit consideration of vehicle dynamics and inter-vehicle interactions, there is no guarantee that the planned arrival times can be upheld in the lower level, where vehicle kinematics and collision avoidance measures are taken into account.
In \citet{Niroumand2020}, an optimization problem is introduced to jointly optimize signal timing, vehicle speeds, and accelerations for mixed-autonomy traffic streams. The formulation becomes tractable due to the choice of the cost function and linearization of nonlinear constraints. It is worth noting that, in this framework, the green duration of each phase must be an integer multiple of the time step size for signal updates. Therefore, the method's performance depends on the level of granularity introduced in signal timing.

Several alternative approaches have been proposed in the literature. Some of these formulations involve optimization problems with non-generalizable cost functions \citep{Yu:2018}, while others rely on heuristics to simplify vehicle trajectory generation \citep{Guo:2019, Feng:2018}. Despite these differences, all these studies have demonstrated that the coordination between traffic lights and vehicles holds promise for improving traffic operations.
\subsection{Contribution of the work}
In this paper, we investigate the optimization-based coordination of traffic lights and vehicle trajectories at a 4-leg intersection, considering various turning options, and we address the shortcomings mentioned in Section \ref{sec:vehiclelight}. Importantly, our approach is adaptable to different intersection layouts. We jointly optimize signal timing and vehicle acceleration, resulting in a mixed-integer nonlinear programming problem.

To tackle this complexity, we propose a two-level decomposition. The upper level addresses green time allocation, while the lower level handles vehicle trajectory optimization problems in parallel. The output of the lower level informs the upper level, and, by employing parametric optimization tools, we solve the upper-level optimization problem.

Within this framework:
\begin{itemize}
\item The cycle length remains flexible, and traffic signal timings are not constrained to integer values, enabling more efficient signal timings.
\item Vehicle trajectories are generated without simplifications or heuristics.
\item There are no constraints on vehicle arrival speeds at the intersection.
\item The general formulation allows for the incorporation of various performance metrics, making it suitable for analyzing intersections with different layouts.
\end{itemize}

This coordinated approach distributes intelligence not only to individual vehicles but also to the traffic lights themselves. It brings together two key players, traffic signals and vehicles, to enhance system performance. Moreover, retaining traffic lights allows vulnerable road users, e.g., pedestrians and cyclists, to be accommodated during the red intervals of conflicting vehicle movements, with the possibility of incorporating additional safety constraints to account for their crossing times and interactions with vehicles.

The remainder of the paper is organized as follows. In Section \ref{sec:problem}, the problem studied in the paper is presented. Section \ref{sec:decomposition} explains the proposed decomposition method followed by the solution approach, which is detailed in Section \ref{sec:solution}. The results of applying the proposed approach are provided in Section \ref{sec:case}, and Section~\ref{sec:conclusion} concludes the paper. 
\section{Problem formulation}\label{sec:problem}
In this section, we treat traffic lights and automated vehicles as agents. We create a model for the intersection scenario, encompassing traffic light and vehicle dynamics, along with constraints for each agent type. The description remains general, with specific details provided in the case study presented in Section \ref{sec:case}. Throughout this article, we work with the following assumptions:

\begin{assumption} 
The control zone, in which the vehicles are controlled, is a circle with a radius of $X$ around the intersection. 
\end{assumption}

\begin{assumption} 
All vehicles are automated and, when in the control zone, they can communicate their location and speed to each other and to the traffic light.
\end{assumption}

\begin{assumption}
The vehicles follow predetermined paths along the road without changing lanes.
\end{assumption}

\begin{assumption} Vehicles making right turns will use the right lane after the intersection, while those making left turns will use the left lane.\label{ass:turn}
\end{assumption}

\begin{assumption} Traffic conditions downstream of the intersection do not restrict vehicle movement, and there is no spill-back effect.
\label{ass:downstream}
\end{assumption}

We note that although assuming fixed vehicle paths and excluding lane changes may be restrictive in general, previous research on signalized intersections has shown that mandatory lane changes can be initiated at a suitable distance upstream of the intersection \citep{Rahimi2024MandatoryLaneChange} and has assumed that such lane changes are completed early, without subsequent discretionary lane changes \citep{E2025TrajectoryOptimisation}. Accordingly, vehicle paths can be treated as fixed within the control zone upstream of the intersection, while no assumptions are made outside this zone.

\subsection{Traffic light timing}
\begin{figure}[h]
\includegraphics{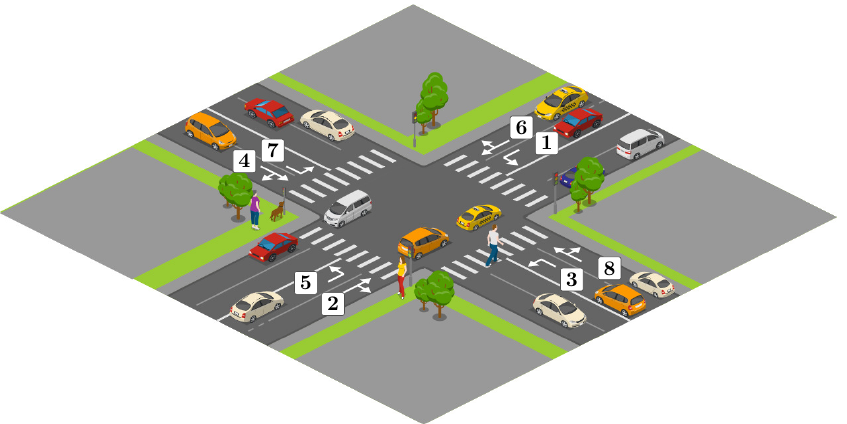}
\caption{Schematic illustration of an intersection. Drawn using https://www.icograms.com.}
\label{fig:intersection}
\end{figure}

We examine a standard 4-leg signalized intersection with two lanes allocated to each of the incoming travel directions, accommodating left-turning and through/right-turning vehicles. Consequently, there are eight approaching movements at the intersection, as illustrated in Figure \ref{fig:intersection}: westbound left (movement 1), eastbound through (movement 2), northbound left (movement 3), southbound through (movement 4), eastbound left (movement 5), westbound through (movement 6), southbound left (movement 7), and northbound through (movement 8).

\begin{figure}[h]
\begin{center}

\includegraphics[scale=0.3]{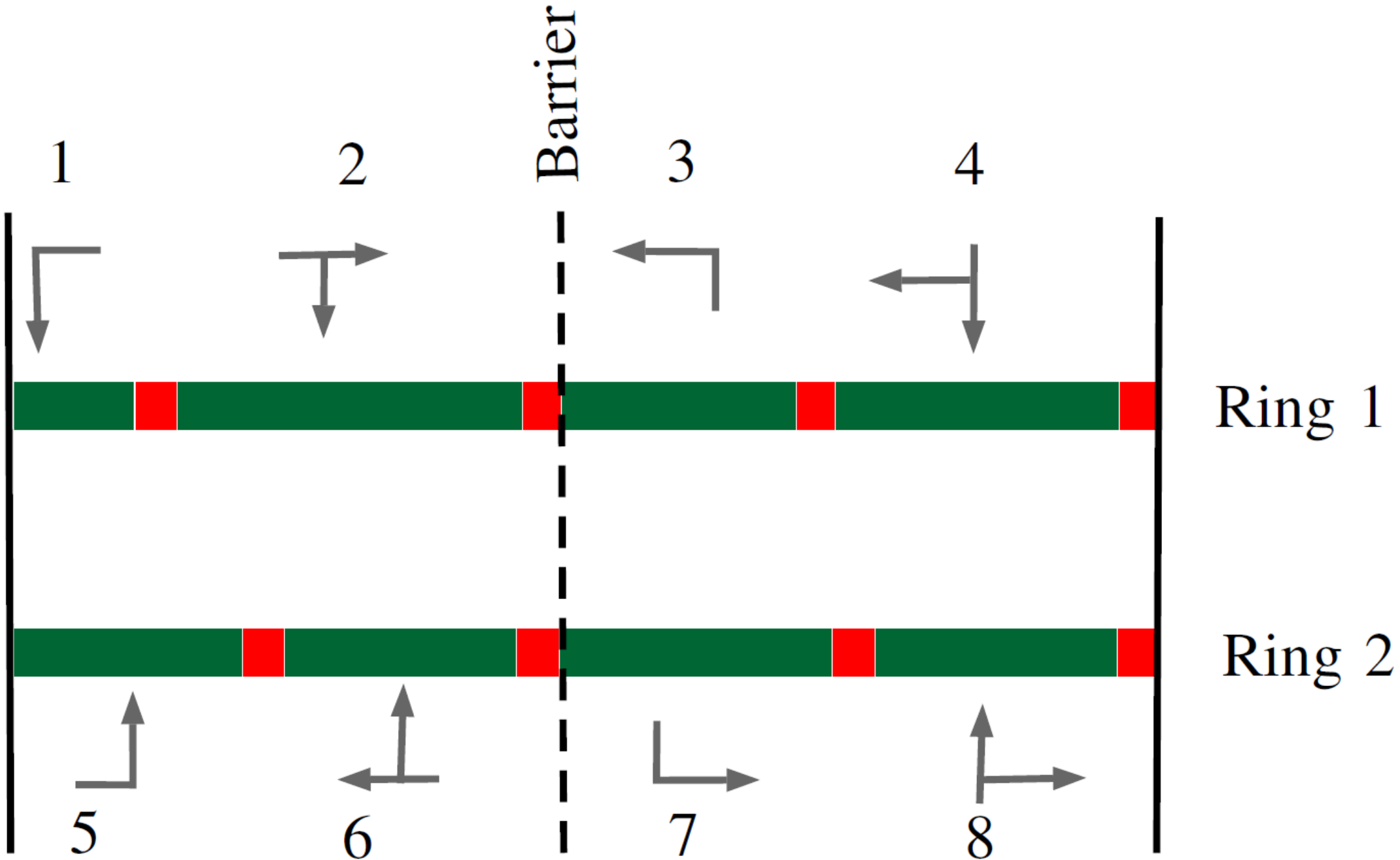}
\caption{Signal phase scheme.}
\label{fig:phase}
\end{center}
\end{figure}

To control the signal timing of each direction, we apply the dual-ring approach shown in Figure \ref{fig:phase}. According to the dual-ring approach, movements 1, 2, 3 and 4 are included in ring 1 and ring 2 includes the remaining movements. The two rings operate independently, except for the timing of the barrier that they share. On either side of the barrier, there are four movements and each of them may have the right-of-way simultaneously with either of the two movements of the other ring on the same side. Thus, the purpose of the traffic signal control is to determine the barrier time and also the green time duration allocated to each of these movements in each cycle. In this work, we assume that no phase is skipped during a traffic light cycle, i.e., all movements are served in each cycle.

Let us denote the times at which the green light for movement $j \in \mathcal{J} = \{\text{1, 2, 3, 4, 5, 6, 7, 8} \}$ in the $c$th cycle starts and ends by $t^{\mathrm{s}_c}_j$ and $t^{\mathrm{e}_c}_j$, respectively. We also collect the start and end times of the green light of movement $j$ in $L$ traffic light cycles with $T_j =(t^{\mathrm{s}_1}_j, \dots, t^{\mathrm{s}_L}_j, t^{\mathrm{e}_1}_j, \dots, t^{\mathrm{e}_L}_j)$ and define the stacked signal timing of all movements in $\mathcal{J} $ as $\bm{T} = (T_{1}, \dots, T_{8})$. Finally, let $\mathcal{F}$ denote the set of all consecutive pairs of movements in each ring; e.g., $\mathcal{F}$ includes the movement pairs (1,2) and (6,7).

According to the dual-ring approach, the allocated timing of all movements must meet certain constraints. These constraints concern the green light duration and some safety measures. To be more specific, the minimum and maximum durations of the green time of each movement in cycle $l$ are limited by $g^\mathrm{{min}}$ and $g^\mathrm{{max}}$, respectively:
\begin{equation}
g^\mathrm{{min}} \leq t^{\mathrm{e}_l}_j - t^{\mathrm{s}_l}_j \leq g^\mathrm{{max}},\quad \forall j\in \mathcal{J}, \forall l\in \mathbb{I}_{[1,L]}, \label{eq:gmin}
\end{equation}
where $ \mathbb{I}_{[1,L]}$ denotes the set of integers in the interval $[1,L]$.  To preserve safety, after the end of the green light of each movement, a clearance time, denoted by $ t^{\mathrm{clr}}$, should be imposed before the next movement in the ring can have the right-of-way: 
\begin{equation}
t^{\mathrm{s}_l}_{j'} - t^{\mathrm{e}_l}_j = t^{\mathrm{clr}},\quad \forall (j,j') \in \mathcal{F}.\label{eq:clearance}
\end{equation}
Moreover, the barrier time in cycle $l$ should be the same in both rings:
\begin{equation}
t_{\mathrm{6}}^{\mathrm{e}_l} = t_{2}^{\mathrm{e}_l} ,\quad \forall l\in \mathbb{I}_{[1,L]}.\label{eq:barrier}
\end{equation}
Similarly, the end of cycle $l$ should be the same in both rings:
\begin{equation}
t_{\mathrm{4}}^{\mathrm{e}_l} = t_{8}^{\mathrm{e}_l} ,\quad \forall l\in \mathbb{I}_{[1,L]}.\label{eq:end}
\end{equation}

\subsection{Vehicle dynamics}

We assume that the motion of vehicle $i$ moving in movement $j$ can be described by
\begin{align}
&\dot{x}_{i,j}^\mathrm{c}(t) = f_{i,j}\bigl(x_{i,j}^\mathrm{c}(t), u_{i,j}^\mathrm{c}(t)\bigr),\label{eq:dynamics-continuous}\\
&g^\mathrm{c}_{i,j}\bigl(x_{i,j}^\mathrm{c}(t), u_{i,j}^\mathrm{c}(t)\bigr) \leq 0,\label{eq:constraint-continuous}
\end{align}
where $x^\mathrm{c}_{i,j}(t)$ and  $u^\mathrm{c}_{i,j}(t) $ denote the state and control input of the vehicle, and the function $g^\mathrm{c}_{i,j}$ captures the possible constraints on the state and the control input. The state vector collects the position, $p^\mathrm{c}_{i,j}(t)$, and speed, $v^\mathrm{c}_{i,j}(t)$, of the vehicle, i.e., $x^\mathrm{c}_{i,j}(t) = \bigl( p^\mathrm{c}_{i,j}(t), v^\mathrm{c}_{i,j}(t) \bigr)$. 

For practical reasons, we formulate the optimization problem in a discrete-time setting. Therefore, with a sampling time $t_\mathrm{s}$, the evolution of the vehicle state and control input during a fixed time horizon $t_\mathrm{f} = M t_\mathrm{s}$ is described by $ x_{i,j} = \bigl(x_{i,j}^{0}, x_{i,j}^1, \dots, x_{i,j}^{M}\bigr)$ and $ u_{i,j} = \bigl(u_{i,j}^{0}, u_{i,j}^{1}, \dots, u_{i,j}^{M}\bigr)$, respectively. For the horizon $t_\mathrm{f}$, the multiple-shooting discretization of  \eqref{eq:dynamics-continuous} and \eqref{eq:constraint-continuous} \citep{NumericalOpt} results in the following equations:
\begin{eqnarray}
&x_{i,j}^{k+1} = F_{i,j}\bigl(x_{i,j}^{k}, u_{i,j}^{k}, t_\mathrm{s}\bigr), \quad &k\in \mathbb{I}_{[0,M-1]}\label{eq:dynamics-dis}\\
&g_{i,j}\bigl(x_{i,j}^{k}, u_{i,j}^{k}\bigr) \leq 0, \quad &k\in \mathbb{I}_{[0,M-1]}\label{eq:constraint-dis}
\end{eqnarray}
where $x_{i,j}^{k} = (p_{i,j}^{k}, v_{i,j}^{k})$. For use in the optimization problem formulation, we define $W_j = (w_{1,j}, \dots, w_{N_j,j})$ to collect the states and control inputs of all $N_j$ vehicles in movement $j$, where $w_{i,j} =( x_{i,j}, u_{i,j})$. We also define $\bm{W} = (W_1, \dots, W_8)$ to collect the states and control inputs of all vehicles in all movements.
\subsection{Collision avoidance}
For all vehicles to move safely in all directions and at all times, two sets of precautionary measures must be taken: one to prevent rear-end collisions between vehicles in the same movement, and one to prevent side collisions between vehicles in conflicting movements within the intersection. In the following, we explain how the former is achieved by keeping a safe headway between vehicles, while the latter is achieved by preventing red-light running. 

To avoid rear-end collisions, each vehicle must keep a safe headway $h$ from its preceding vehicle. Hence, if vehicle $i$ moves in front of vehicle $i+1$ in movement $j$, the following inequality constraint must hold for their positions:
\begin{equation}\label{eq:rear-end}
p_{i+1,j}^{k} +h \leq p_{i,j}^{k},\quad  \forall i\in\mathbb{I}_{[1,N_{j}-1]},  \forall k\in \mathbb{I}_{[0,M]}.
\end{equation}
Note that, to extend the formulation to cases where downstream vehicle breakdowns or unexpected delays invalidate Assumption \ref{ass:downstream}, collision avoidance constraints analogous to \eqref{eq:rear-end} can also be imposed on vehicles entering the same downstream lane. The positions of these vehicles would be expressed along a common downstream lane, and the safe headway $h$ would be imposed between preceding and following vehicles regardless of their originating movements.

Vehicles moving toward the intersection should only pass the intersection when they have the right-of-way. Let $p_{i,j}(\tau,w_{i,j})$ denote the position of vehicle $i$ in movement $j$ at time $\tau$, defined as 
\begin{equation}
p_{i,j}(\tau,w_{i,j}) \coloneqq [1, 0] \, F_{i,j}\bigl(x_{i,j}^{k}, u_{i,j}^{k}, \tau - kt_\mathrm{s}\bigr),
\end{equation}
where $k = \floor{\frac{\tau}{t_\mathrm{s}}}$ and $\floor{.}$ denotes the floor function. Note that $\tau$ need not be an integer multiple of $t_\mathrm{s}$. A necessary set of conditions that allows $r_j$ out of $N_j$ vehicles in movement $j$ to pass the intersection during the green light period in the $c$th cycle is
\begin{align}
p_{i,j}(t^{\mathrm{s}_c}_j, w_{i,j})& \leq X, \quad \forall i \in \mathbb{I}_{[1,N_j]},\label{eq:green1}\\
p_{i,j}(t^{\mathrm{e}_c}_j, w_{i,j}) &\geq X, \quad\forall i \in \mathbb{I}_{[1,r_j]}, \label{eq:green2}\\
p_{i,j}(t^{\mathrm{s}_{c+1}}_j, w_{i,j}) &\leq X, \quad\forall i \in \mathbb{I}_{[r_j+1, N_j]}. \label{eq:green3}
\end{align}
Equation \eqref{eq:green1} ensures that no vehicle in movement $j$ crosses the intersection before the corresponding light turns green. By meeting \eqref{eq:green2}, the first $r_j$ vehicles in movement $j$ pass the intersection before the light turns red. Complying with \eqref{eq:green3} prevents the remaining vehicles from crossing the intersection until the start of movement $j$'s green light period in the next cycle.
\subsection{Coordination problem}
With all these ingredients in place, the optimal control formulation of the coordination problem between vehicles and traffic lights reads:
\begin{subequations}\label{eq:min1}
\begin{align}
\min_{\bm{T},\bm{W},\bm{r}} \quad &O(\bm{W})\\
\text{s.t.} \quad&\eqref{eq:gmin},\eqref{eq:clearance}, \eqref{eq:barrier}, \eqref{eq:end},\\
&\eqref{eq:dynamics-dis}, \eqref{eq:constraint-dis},\quad \forall j \in \mathcal{J}, \forall i\in\mathbb{I}_{[1,N_j]},\\
& \eqref{eq:rear-end}, \eqref{eq:green1}, \eqref{eq:green2},\eqref{eq:green3},\quad  \forall j \in \mathcal{J},
\end{align}
\end{subequations}
where $\bm{r} =(r_1,\dots, r_8)$ and $O(\bm{W}) = \sum_{j=1}^{8} \sum_{i=1}^{N_j} o_{i,j}(w_{i,j})$, in which $o_{i,j}(w_{i,j})$ denotes the objective function of vehicle $i$ in movement $j$, defined as
\begin{equation}
o_{i,j}(w_{i,j}) =  \sum_{k=0}^{M-1}l_{i,j}(x_{i,j}^k, u_{i,j}^k).
\end{equation}
Here, the  twice continuously differentiable function $l_{i,j}(\cdot,\cdot)$ describes the performance
metric of the vehicle.

Problem \eqref{eq:min1} is a challenging mixed-integer nonlinear program, especially for large instances. To address this challenge, we outline a decomposition method for obtaining an approximate solution to problem \eqref{eq:min1}.

\begin{remark}
To extend the proposed framework to mixed traffic with human-driven and automated vehicles, prediction of human-driven vehicle trajectories would need to be integrated into the formulation. Inspiration may be taken from \citet{Kojchev2025MixedTraffic} and \citet{Chen2021MixedPlatoon}. The resulting predicted trajectories can then be incorporated into the coordination problem, while optimizing the traffic light timing and the control inputs of automated vehicles. Appropriate safety margins in the collision avoidance constraints would also be needed to account for prediction uncertainty.
\end{remark}
\section{Decomposition}\label{sec:decomposition}
In our two-level decomposition method, the upper level addresses the green time allocation problem for all movements, given $\bm{W}$ and $\bm{r}$. The lower level involves an optimization problem for each movement. Specifically, problem \eqref{eq:min1} can be decomposed as follows:

\begin{subequations}\label{eq:min_decompose1}
\begin{align}
\min_{\bm{T}} \quad &\sum_{j=1}^{8}\Phi_j(T_j)\\
\text{s.t.} \quad&\eqref{eq:gmin},\eqref{eq:clearance}, \eqref{eq:barrier}, \eqref{eq:end},
\end{align}
\end{subequations}
where $\Phi_j(T_j)$ is the optimal value of the following parametric optimization problem for the $j$th movement:
\begin{subequations}\label{eq:min_decompose2}
\begin{align}
\Phi_j(T_j) = \min_{r_j, W_j} \quad &\sum_{i=1}^{N_j}o_{i,j}(w_{i,j})\label{eq:min_decompose2_1}\\
\text{s.t.} \quad&\eqref{eq:dynamics-dis}, \eqref{eq:constraint-dis},\quad  \forall i\in\mathbb{I}_{[1,N_j]},\\
& \eqref{eq:rear-end}, \eqref{eq:green1}, \eqref{eq:green2},\eqref{eq:green3}.  
\end{align}
\end{subequations}
Note that, in general, \eqref{eq:min_decompose2} remains a mixed-integer nonlinear programming (NLP) problem due to the integer nature of $r_j$, representing the number of vehicles passing the current or the first upcoming green light in movement $j$. However, $\Phi_j(T_j) $ can be found equivalently by solving $N_j+1$ independent NLP problems in parallel, corresponding to all $N_j+1$ possible values that $r_j$ may take, and selecting the one with the minimum cost. That is to say,
\begin{equation}\label{eq:minfunc}
\Phi_j(T_j) = \min\left\{\tilde{\Phi}_{j_0}(T_j),\dots,\tilde{ \Phi}_{j_{N_j}}(T_j)\right\}.
\end{equation}
The cost $\tilde{\Phi}_{j_b}(T_j)$, for $b\in \mathbb{I}_{[0,N_j]}$, is the solution of the optimization problem \eqref{eq:min_decompose2} when $r_j =b$,  i.e., when $b$ vehicles in the $j$th movement pass the intersection during the corresponding green interval of the $c$th cycle. This optimization problem reads 
\begin{subequations}\label{eq:min_decompose3}
\begin{align}
\tilde{\Phi}_{j_b}(T_j) =\nonumber\\
 \min_{ W_j} \quad &\sum_{i=1}^{N_j}o_{i,j}(w_{i,j})\\
\text{s.t.} \quad&\eqref{eq:dynamics-dis}, \eqref{eq:constraint-dis},\quad  \forall i\in\mathbb{I}_{[1,N_j]},\\
& \eqref{eq:rear-end}, \\
&p_{i,j}(t^{\mathrm{s}_c}_j, w_{i,j}) \leq X_j, \quad \forall i \in \mathbb{I}_{[1,N_j]},\label{eq:green1_b}\\
&p_{i,j}(t^{\mathrm{e}_c}_j, w_{i,j}) \geq X_j, \quad\forall i \in \mathbb{I}_{[1,b]}, \label{eq:green2_b}\\
&p_{i,j}(t^{\mathrm{s}_{c+1}}_j, w_{i,j}) \leq X_j, \quad\forall i \in \mathbb{I}_{[b+1, N_j]}. \label{eq:green3_b}  
\end{align}
\end{subequations}

Having the solution of problem \eqref{eq:min_decompose3} for all values of $b$, $\Phi_j(\cdot)$ can be obtained using \eqref{eq:minfunc}. The resulting function is non-smooth and hence not suitable for the gradient-based algorithm that is used to solve \eqref{eq:min_decompose1}. Therefore, \eqref{eq:minfunc} is approximated by the following smooth function:
\begin{equation}\label{eq:logapproximation}
\Phi_j(T_j) \approx -\frac{1}{\beta} \log\sum_{b=0}^{N_j}\exp\left(-\beta \,\tilde{\Phi}_{j_b}(T_j)\right),
\end{equation}
for $\beta>1$ \citep{Jamshidnejad:2018aa}.

Let us summarize the proposed decomposition method. To provide the cost function of the optimization problem in the upper level, $\sum_{j=1}^{8}(N_j+1)$ optimization problems of type \eqref{eq:min_decompose3} must be solved in the lower level. Intuitively, each movement needs to solve its corresponding optimization problem to inform the upper-level controller of the minimum trajectory cost of vehicles approaching the intersection, given the current traffic light timing and the number and locations of vehicles within that movement's control zone. Note that the underlying optimization problems of the movements are fully independent of each other.

Namely, for each movement $j$, given the initial value of $W_j$ from all vehicles in the movement, the optimization problem \eqref{eq:min_decompose3} should be solved for all possible values of $b$ in order to calculate $\Phi_j(T_j)$ using \eqref{eq:logapproximation}. Note that all of these optimization problems can be solved in parallel. We defer the discussion concerning the computational demand of the approach to Section \ref{sec:computation}.

Having $\Phi_j(T_j)$ for all $j$, the following section explains how tools from parametric optimization can be used to solve \eqref{eq:min_decompose1}.

\section{Solution approach}\label{sec:solution}
In this section, we explain a numerical algorithm to solve problem \eqref{eq:min_decompose1} by using sequential quadratic programming (SQP). SQP is selected because of its faster convergence compared with first-order methods and because the derivatives that it requires can be obtained using tools from parametric optimization. For notational convenience, let us rewrite problem \eqref{eq:min_decompose1} in a more compact form as
\begin{subequations}\label{eq:condense}
\begin{align}
\min_{\bm{T}} \quad &\sum_{j=1}^{8}\Phi_j(T_j),\\
\text{s.t.} \quad&C\bm{T}\leq 0,\label{eq:condense_ineq}\\
&D\bm{T} =0,\label{eq:condense_eq}
\end{align}
\end{subequations}
where \eqref{eq:condense_ineq} and \eqref{eq:condense_eq} respectively collect all the inequality and equality constraints associated with \eqref{eq:gmin}--\eqref{eq:end}. Using the SQP algorithm, the solution to \eqref{eq:condense} is obtained through the following iterations:
\begin{align}
\bm{T}^{+} &=\bm{T} + \alpha \Delta \bm{T},\\
{\mu}^{+} &= (1-\alpha){\mu} + \alpha \mu^{\mathrm{QP}},\\
{\lambda}^{+} &=(1-\alpha){\lambda} + \alpha \lambda^{\mathrm{QP}},
\end{align}
where $(.)^{+}$ denotes the value of the corresponding variable in the next iteration, $\mu$ and $\lambda$ are the vectors of multipliers associated with constraints \eqref{eq:condense_ineq} and \eqref{eq:condense_eq}, respectively, $\alpha \in (0,1]$ is the step size, which can be adjusted using available algorithms such as line search, and $\Delta \bm{T}$, $\mu^{\mathrm{QP}}$, and $\lambda^{\mathrm{QP}}$ are the primal-dual solution of the following quadratic programming (QP) problem:
\begin{subequations}\label{eq:QP}
\begin{align}
\min_{\Delta\bm{ T}} \quad &\frac{1}{2}  \Delta\bm{ T}^{\top}H(\bm{T},\mu, \lambda) \Delta\bm{ T}  + \bm{\Phi}(\bm{T})\Delta\bm{T}  \label{eq:QP1}     \\
\text{s.t.} \quad&C\bm{T} + C\Delta\bm{T}\leq 0,\label{eq:QP2}\\
&D\bm{T}+ D\Delta\bm{T} =0,\label{eq:QP3}
\end{align}
\end{subequations}
with 
\begin{equation}\label{eq:gradient}
\bm{\Phi}(\bm{T}) = [\nabla^{\top}\Phi_1(T_1), \nabla^{\top}\Phi_2(T_2),\dots, \nabla^{\top}\Phi_8(T_8)].
\end{equation}
In \eqref{eq:QP1}, $H$ is the Hessian of the Lagrange function of problem \eqref{eq:condense}, which is defined as
\begin{equation}
\mathcal{L} (\bm{T},\mu, \lambda)= \sum_{j=1}^{8}\Phi_j(T_j) + \mu^{\top}C\bm{T} + \lambda^{\top}D\bm{T}.
\end{equation}
Note that, since all constraints of problem \eqref{eq:condense} are linear and its objective function is separable in $T_1,\dots,T_8$, the Hessian $H$ is a block diagonal matrix with its $j$th block obtained from
\begin{equation}\label{eq:hessian}
H_j = \nabla^2\Phi_j(T_j).
\end{equation}
where regularization is used when $H_j$ is not positive definite \citep{NumericalOpt}.

From \eqref{eq:gradient} and \eqref{eq:hessian}, it is clear that the first-order and second-order derivatives of $\Phi_j(T_j)$ are required for solving the QP problem \eqref{eq:QP} and, subsequently, for finding the solution of \eqref{eq:condense}. The following subsection details the required steps for computing these derivatives.
\subsection{Computing derivatives}

If the $o$th entry of $T_j$ is denoted by $T^o_j$, then, from \eqref{eq:logapproximation} and as derived in Appendix \ref{app:sumlog}, the first-order and second-order derivatives of  $\Phi_j(T_j)$ with respect to $T^o_j$ are given by
\begin{equation}\label{eq:exp_first_derivative}
\dv{\Phi_j(T_j)}{T^o_j} =  \displaystyle\frac{\displaystyle \sum_{b=0}^{N_j} \displaystyle   \dv{\tilde{\Phi}_{j_b}(T_{j})}{T^o_j}  \exp\left(-\beta\tilde{\Phi}_{j_b}(T_{j})\right)}{{\displaystyle\sum_{b=0}^{N_j}  \exp\left(-\beta\tilde{\Phi}_{j_b}(T_{j})\right)}},
\end{equation}
and
\begin{align}\label{eq:exp_second_derivative}
&\dv[2]{\Phi_j(T_j)}{{T^o_j}} = \displaystyle \beta \left( \frac{\displaystyle \sum_{b=0}^{N_j} \displaystyle \dv{\tilde{\Phi}_{j_b}(T_{j})}{T^o_j} \exp \left(-\beta\tilde{\Phi}_{j_b}(T_{j})\right)} {\displaystyle\sum_{b=0}^{N_j}  \exp\left(-\beta\tilde{\Phi}_{j_b}(T_{j})\right)}   \right)^2 \nonumber\\
 & + \frac{ \displaystyle\! \sum_{b=0}^{N_j}\! \exp \left(-\beta\tilde{\Phi}_{j_b}(T_{j})\right) \!\left(\! \dv[2]{\tilde{\Phi}_{j_b}(T_{j})}{{T^o_j}} -\beta \left(\dv{\tilde{\Phi}_{j_b}(T_{j})}{T^o_j}\right)^2\!\right)}             {\displaystyle\sum_{b=0}^{N_j} \! \exp \left(-\beta\tilde{\Phi}_{j_b}(T_{j})\right)}.
\end{align}
Therefore, in order to calculate the terms in \eqref{eq:exp_first_derivative} and \eqref{eq:exp_second_derivative}, the first-order and second-order derivatives of the solution of \eqref{eq:min_decompose3} with respect to $T^o_j$ (and in general $T_j$) are required.  We note that, since $o_{i,j}$ and the constraints in \eqref{eq:min_decompose3} are twice differentiable, and assuming that, in \eqref{eq:min_decompose3}, at the solution $\bar{W}_j(T_j)$, the linear independence constraint qualification (LICQ) and the second-order sufficient condition (SOSC) \citep{NumericalOpt} hold $\forall T_{j}\in \mathcal{T}$, $\tilde{\Phi}_{j_b}(\cdot)$ will be piecewise twice continuously differentiable on $\mathcal{T}$. At the points where an inequality constraint is weakly active, sensitivities may be defined only up to the subgradients generated by the possible active sets \citep{sensitivity}. To address this in practice, interior-point techniques can be employed when solving \eqref{eq:min_decompose3}.

\subsubsection{First-order derivatives}For the sake of explanation, let us choose $T^o_j $ to be $t^{\mathrm{s}_c}_j$. If $\mathcal{L}_{j_b}$ denotes the Lagrange function of \eqref{eq:min_decompose3}, $\bar{W}_j$ is the primal solution of \eqref{eq:min_decompose3} for a given $T_j$,  $P_j(t^{\mathrm{s}_c}_j, \bar{W}_j)$ is the stacked vector of the constraints \eqref{eq:green1_b} for all $i$, and $\gamma$ is the vector of Lagrange multipliers associated with these constraints, then the first-order derivative of $\tilde{\Phi}_{j_b}(T_{j})$, for any $b\in \mathbb{I}_{[0,N_j]}$, with respect to $t^{\mathrm{s}_c}_j$ is given by \citep{sensitivity}:

\begin{equation}\label{eq:first_derivative}
\dv{\tilde{\Phi}_{j_b}(T_{j})}{t^{\mathrm{s}_c}_j} = \pdv{\mathcal{L}_{j_b}}{t^{\mathrm{s}_c}_j} = \gamma^{\top} \pdv{P_j(t^{\mathrm{s}_c}_j, \bar{W}_j)}{t^{\mathrm{s}_c}_j}.
\end{equation}
\subsubsection{Second-order derivatives} Differentiating \eqref{eq:first_derivative} results in
 
\begin{align}\label{eq:second_derivative}
 \dv[2]{\tilde{\Phi}_{j_b}}{{t^{\mathrm{s}_c}_j}} = \left(  \dv{\gamma}{t^{\mathrm{s}_c}_j}\right)^{\top}     \pdv{P_j}{t^{\mathrm{s}_c}_j} +  
 \gamma^{\top} \left( \pdv[2]{{P_j}}{{t^{\mathrm{s}_c}_j}}   + \pdv{P_j}{t^{\mathrm{s}_c}_j}{\bar{W}_j}  \dv
{\bar{W}_j}{t^{\mathrm{s}_c}_j}   \right)
\end{align} 
where the arguments in \eqref{eq:second_derivative} have been dropped for brevity. 

As seen in \eqref{eq:second_derivative}, calculating the second-order derivative of $\tilde{\Phi}_{j_b}$  hinges on finding the sensitivity of the primal solution ${\bar{W}_j}$ and the Lagrange multiplier ${\gamma}$ with respect to the corresponding parameter. To obtain the value of these terms, we make use of the Karush--Kuhn--Tucker (KKT) matrix of problem \eqref{eq:min_decompose3} and the KKT optimality conditions. The optimality condition is briefly explained in Appendix \ref{app:KKT}, and interested readers are referred to \citet{NumericalOpt} for more details.  

If the primal-dual solution of \eqref{eq:min_decompose3} for a given $T_j$ is denoted by $\bar{Y}(T_j) $, the KKT optimality condition at this point can be compactly written as:
\begin{equation}\label{eq:KKT}
K\left(\bar{Y}(T_j), T_j\right) = 0.
\end{equation} 
Thus, differentiating $K$ yields a system of linear equations for the sensitivity differentials of the primal and dual solutions. For example, to find the value of  $\dv{\bar{W}_j}{t^{\mathrm{s}_c}_j}$ and $\dv{\gamma}{t^{\mathrm{s}_c}_j}$, which are embedded in $\dv{Y}{t^{\mathrm{s}_c}_j}$, it is sufficient to take the derivative of \eqref{eq:KKT} with respect to $t^{\mathrm{s}_c}_j$ to obtain
\begin{equation}\label{eq:KKT-derivative}
\dv{K}{t^{\mathrm{s}_c}_j} = \pdv{K}{t^{\mathrm{s}_c}_j} + \left(\nabla_Y K\right)^{\top} \dv{Y}{t^{\mathrm{s}_c}_j} = 0,
\end{equation} 
and solve this system of linear equations to find $\dv{Y}{t^{\mathrm{s}_c}_j}$. Note that if \eqref{eq:min_decompose3} is solved using a second-order method, as in this work, the solver needs to factorize the KKT Jacobian to perform its iterations. Specifically, following the notation of \eqref{eq:KKT}, the solver takes steps of the form
\begin{equation*}
Y^{c+1} = Y^c - \alpha^c
\left(\frac{\partial K}{\partial Y}(Y^c,T_j)\right)^{-1} K(Y^c,T_j),
\end{equation*}
where $c$ denotes the iteration index and $\alpha^c$ is the step size. This indicates that the solver factorizes the matrix $\partial K/\partial Y$ during its iterations. When the factorization at the solution is available from the solver, it can therefore be reused to solve the sensitivity system \eqref{eq:KKT-derivative}, reducing the additional computational cost of evaluating the derivatives in \eqref{eq:first_derivative} and \eqref{eq:second_derivative}.

\section{Evaluation scenarios and results}\label{sec:case}
This section assesses the performance of the proposed approach through an illustrative case study and compares it to a baseline approach. We begin by introducing the motion model, the primary and secondary control objectives, and the performance scores used in the case studies. The parameters used in the case studies are listed in Table~\ref{table:parameters}.

\subsection{Motion model}

All controllers use the double-integrator longitudinal dynamics, i.e., $f_{i,j}\bigl(x^\mathrm{c}_{i,j}(t),u^\mathrm{c}_{i,j}(t)\bigr)$ in \eqref{eq:dynamics-continuous} is defined as follows:
\begin{equation}
f_{i,j}\bigl(x^\mathrm{c}_{i,j}(t),u^\mathrm{c}_{i,j}(t)\bigr) = 
\begin{bmatrix}
\dot{p}^\mathrm{c}_{i,j}(t)\\
\dot{v}^\mathrm{c}_{i,j}(t)\\
\end{bmatrix}
=
\begin{bmatrix}
v^\mathrm{c}_{i,j}(t)\\
u^\mathrm{c}_{i,j}(t)\\
\end{bmatrix}.
\end{equation}
Although the prediction model is simple, experimental results indicate that it is sufficient for similar applications \citep{Hult:2019aa}.
\begin{table*}
\centering
\caption{Simulation parameters.}
\label{table:parameters}
\small
\setlength{\tabcolsep}{0pt} 

\begin{tabular*}{\textwidth}{@{\extracolsep{\fill}} ccccccccccc @{}}
\toprule
$X$ & $v^\mathrm{max}$ & $v^\mathrm{min}$ & $v^\mathrm{d}$ & $u^\mathrm{min}$ & $u^\mathrm{max}$ & $t_\mathrm{s}$ & $h$ & $t^{\mathrm{clr}}$ & $g^\mathrm{max}$ & $g^\mathrm{min}$ \\
(\unit{m}) & (\unit{m/s}) & (\unit{m/s}) & (\unit{m/s}) & (\unit{m/s^2}) & (\unit{m/s^2}) & (\unit{s}) & (\unit{m}) & (\unit{s}) & (\unit{s}) & (\unit{s}) \\
\midrule
200 & 14 & 0 & 12.27 & -2 & 2 & 2 & 3 & 2 & 25 & 6 \\
\bottomrule
\end{tabular*}

\vspace{0.6cm} 

\begin{tabular*}{\textwidth}{@{\extracolsep{\fill}} cccccccccccc @{}}
\toprule
$Q^\mathrm{u}$ & $Q^\mathrm{v}$ & $m$ & $A_\mathrm{f}$ & $C_\mathrm{D}$ & $\rho$ & $C_\mathrm{r}$ & $c_1$ & $c_2$ & $\eta$ & $\alpha_1$ & $\alpha_2$ \\
-- & -- & (\unit{kg}) & (\unit{m^2}) & -- & (\unit{kg/m^3}) & -- & -- & -- & -- & -- & -- \\
\midrule
2.5 & 0.0153 & 1487 & 2.3 & 0.3 & 1.226 & 1.75 & 0.0328 & 4.575 & 0.92 & 4.29e-5 & 1e-6 \\
\bottomrule
\end{tabular*}
\end{table*}
\subsection{Main objective function}
We aim to address the coordination of traffic lights and vehicle speeds, focusing on the conventional efficiency metric of acceleration variation while promoting comfortable cruising at a predefined desired speed.

The corresponding objective function $o_{i,j}(w_{i,j})$ in \eqref{eq:min_decompose2_1} is therefore taken to be
\begin{equation}\label{eq:cost_fuel}
o_{i,j}(w_{i,j}) = \sum_{k=0}^{M-1} \Bigl(Q^\mathrm{u}(u_{i,j}^k)^2 + Q^\mathrm{v} (v_{i,j}^k - v^\mathrm{d}_{i,j})^2\Bigr),
\end{equation}
Here, $Q^\mathrm{u}$ and $Q^\mathrm{v}$ are objective weights, and ${v^\mathrm{d}_{i,j}}$ represents the known desired speed of the $i$th vehicle in the $j$th movement. The second term in \eqref{eq:cost_fuel} accounts for individual vehicle preferences during trajectory optimization. The corresponding performance score for the quadratic objective \eqref{eq:cost_fuel} is calculated as the average value across all vehicles during the simulation.

\begin{equation}\label{eq:score_performance}
S^\mathrm{c} = \frac{1}{|\mathcal{N}|} \sum_{(i,j)\in \mathcal{N}}  \sum_{k=k^\mathrm{e}_{i,j}}^{k^\mathrm{l}_{i,j}}\Bigl(Q^\mathrm{u}(u_{i,j}^k)^2 + Q^\mathrm{v} (v_{i,j}^k - v^\mathrm{d}_{i,j})^2\Bigr),
\end{equation}
where $\mathcal{N}$ collects the indices of all the vehicles that entered and left the control zone during the simulation time, and $k^\mathrm{e}_{i,j}$ and $k^\mathrm{l}_{i,j}$ denote the time steps at which vehicle $i$ in movement $j$ enters and leaves the control zone, respectively.
\subsection{Secondary objective functions}
In all case studies, the objective function used in the optimization problem follows \eqref{eq:cost_fuel}. Additionally, we use total travel time and precise fuel consumption estimation as performance metrics to assess and compare the proposed control algorithm.
\subsubsection{Total travel time} We measure the total travel time of a vehicle $i$ in movement $j$ as the total time it takes the vehicle to travel through the control zone. If we denote this time by $\bar{t}_{i,j}$, the corresponding performance score reads 
\begin{equation}\label{eq:score_time}
S^\mathrm{t} = \frac{1}{|\mathcal{N}|} \sum_{(i,j)\in \mathcal{N}} \bar{t}_{i,j}.
\end{equation}
Although not explicitly, the cost function in \eqref{eq:cost_fuel} indirectly incorporates travel time minimization through its second term, which penalizes low velocities when the desired velocity $v^\mathrm{d}_{i,j}$ takes a relatively high value.

\begin{table}
\small
\caption{Demand  (\si{veh\per\hour\per lane}) in each movement corresponding to the Demand Profiles I, II, and III. }
\begin{center}
\begin{tabularx}{\columnwidth}{Xcccccccc}

\rot[90]{Demand Profile} &\rot[90]{Movement 1}&\rot[90]{Movement 2}& \rot[90]{Movement 3}& \rot[90]{Movement 4}&\rot[90]{Movement 5}& \rot[90]{Movement 6}& \rot[90]{Movement 7}& \rot[90]{Movement 8}\\[.02cm]\hline\hline
\rule{0pt}{1.1\normalbaselineskip} I& 800&800  &800&800&800&800&800&800 \\[.2cm]\hline
\rule{0pt}{1\normalbaselineskip} II&1000 &1000  &1000  &1000 & 1000& 1000 &1000&1000 \\
\rule{0pt}{1\normalbaselineskip} III&900 &1500  &400  &400 & 900& 1500 &400&400 \\
\end{tabularx}
\end{center}
\label{table:demand}
\end{table}
\begin{table}
\small
\caption{Simulation cases and their specifications.}
\begin{center}
\begin{tabularx}{\columnwidth}{c>{\setlength\hsize{1\hsize}\centering}Xcc}
& The proposed control scheme & Demand Profile & Cycle length\\\hline\hline
\rule{0pt}{1.1\normalbaselineskip}Case 1& yes& I & not-fixed\\
\rule{0pt}{1.1\normalbaselineskip}Case 2& no & I & 54\\
\rule{0pt}{1.1\normalbaselineskip}Case 3& no&I&60\\
\rule{0pt}{1.1\normalbaselineskip}Case 4& yes &II& not-fixed\\
\rule{0pt}{1.1\normalbaselineskip}Case 5& no & II & 54\\
\rule{0pt}{1.1\normalbaselineskip}Case 6& no & II & 60\\
\rule{0pt}{1.1\normalbaselineskip}Case 7& yes &III& not-fixed\\
\rule{0pt}{1.1\normalbaselineskip}Case 8&no&III&54\\
\rule{0pt}{1.1\normalbaselineskip}Case 9&no&III&60\\
\end{tabularx}
\end{center}
\label{table:scenarios}
\end{table}

\subsubsection{Explicit fuel consumption }
For a more accurate estimation of total fuel consumption, we employ the Virginia Tech comprehensive power-based fuel consumption model (VT-CPFM), a second-order polynomial model of the vehicle's power~\citep{VirginiaTech}. According to VT-CPFM, if the power of vehicle $i$ in movement $j$ at time $t$ is represented as $\mathcal{P}_{i,j}(t)$ (kW), the fuel consumption rate of the vehicle, denoted by $\mathcal{F}_{i,j}(t)$, can be expressed as follows:
\begin{equation}\label{eq:VT-fuel}
\mathcal{F}_{i,j}(t) =
\begin{cases}
 \alpha_0 + \alpha_1 \mathcal{P}_{i,j}(t) + \alpha_2 \bigl(\mathcal{P}_{i,j}(t)\bigr)^2 & \forall \; \mathcal{P}_{i,j}(t) \geq 0\\\
\alpha_0 &  \forall \; \mathcal{P}_{i,j}(t) < 0
\end{cases}
\end{equation}
where $\alpha_0$, $\alpha_1$ and $\alpha_2$ are vehicle-specific model constants. The power exerted at any instant $t$ can be computed as
\begin{equation}\label{eq:VT-power}
\mathcal{P}_{i,j}(t) = \frac{R_{i,j}(t) + 1.04\, m u_{i,j}(t)}{1000 \eta}v_{i,j}(t),
\end{equation}
where $m$ is the vehicle mass (\si{kg}), and $\eta$ is the driveline efficiency. The resistance force $R_{i,j}(t)$ in \eqref{eq:VT-power} is computed as the sum of the aerodynamic, rolling, and grade resistance forces and reads
\begin{equation}\label{eq:VT-R}
R_{i,j}(t) = 0.5\, \rho\, C_\mathrm{D}A_\mathrm{f}\bigl(v_{i,j}(t)\bigr)^2 +  \frac{C_\mathrm{r}}{1000}m\,g\bigl(c_1v_{i,j}(t) + c_2\bigr) +  mg\,\tan(\theta).
\end{equation}
Here, $\rho$ represents air density (\si{kg/m^3}), $C_\mathrm{D}$ is the vehicle drag coefficient (unitless), $A_\mathrm{f}$ is the vehicle frontal area (\si{m^2}), and $g$ is the gravitational acceleration (9.8066 \si{m/s^2}). Additionally, $C_\mathrm{r}$, $c_1$, and $c_2$ are rolling resistance parameters. It is worth noting that the vehicle parameters in \eqref{eq:VT-power} and \eqref{eq:VT-R} are vehicle-specific, but for clarity, we omit their indices in these equations. Furthermore, $\theta$ in \eqref{eq:VT-R} represents the road grade angle, which is assumed to be zero. The corresponding performance score for assessing the fuel consumption of all vehicles during the simulation time is defined as:
\begin{equation}\label{eq:score_fuel}
S^\mathrm{e} = \frac{1}{|\mathcal{N}|} \sum_{(i,j)\in \mathcal{N}}  \sum_{k=k^\mathrm{e}_{i,j}}^{k^\mathrm{l}_{i,j}} {\mathcal{F}}_{i,j}^k,
\end{equation}
where ${\mathcal{F}}_{i,j}^k$ is the fuel consumption rate of vehicle $i$ in movement $j$ at time step $k$.

\subsection{Computational demand}\label{sec:computation}
Evaluating the cost function of the upper-level problem \eqref{eq:min_decompose1} requires solving $\sum_{j=1}^{8}(N_j+1)$ optimal control problems (OCPs) of the form \eqref{eq:min_decompose3}. An increase in the number of vehicles increases both the number of OCPs and the number of state and control variables in each OCP. Since the OCPs for different movements $j$ and values of $b$ are independent, they can be solved in parallel. In the following, we explain how parallel computation, offline preparation, and the use of parametric optimization can be combined to reduce the total computation time.
\subsubsection{Parallel computation}
Most of the optimal control problems formulated in Section \ref{sec:decomposition} can be solved in parallel. Specifically, the OCPs in \eqref{eq:min_decompose3}, formulated for all movements $j$ and all possible values of $b$ within each movement, are independent of each other and can be solved concurrently. Let $\tilde{t}^j_b$ denote the computation time for solving \eqref{eq:min_decompose3} for the pair $(j,b)$. With sufficient computational resources for fully parallel execution, and neglecting the overhead associated with parallel computation, the time required to solve these OCPs for one evaluation of the cost function in \eqref{eq:min_decompose1} is
\begin{equation}\label{eq:parallel-time}
 \max_{j\in\mathcal{J},b\in\mathbb{I}_{[0,N_j]}}\tilde{t}^j_b.
\end{equation}
Thus, the time required to solve all OCPs of the form \eqref{eq:min_decompose3} is reduced from the sum of their individual solve times under sequential computation to the maximum in \eqref{eq:parallel-time} under fully parallel computation. Parallel computation is particularly useful under high demand conditions, since the additional OCPs associated with larger values of $N_j$ remain independent.

We also note that, although in theory, $b$ in \eqref{eq:parallel-time} is selected from $\{0,\ldots,N_j\}$, simple heuristics can be used in practice to reduce the number of optimization problems that need to be solved. More specifically, recall that $b$ for each movement represents the number of vehicles in that movement that will pass during the current green phase. In each movement, given the state of the vehicles and considering the maximum acceleration and deceleration, two categories of vehicles can be readily identified: 1) vehicles that are too close to the intersection to stop safely before the stop line, even under maximum deceleration, and therefore must pass during the current green phase and 2) vehicles that are too far from the intersection to reach it during the current green phase, even under maximum acceleration, and therefore must wait for the next green phase. For example, suppose that, at a sampling time, a movement contains 12 vehicles, of which 4 belong to the first category and 4 belong to the second. Then, it is sufficient to solve only 5 optimization problems, corresponding to $b\in\{4,\ldots,8\}$, rather than considering all 13 possible values of $b\in\{0,\ldots,12\}$.

\subsubsection{Offline preparation phase}

The optimization problem in \eqref{eq:min_decompose3} is essentially a parametric NLP. The parameters include the initial value of the optimal control problem, i.e., the current speed and position of vehicles and the timing of the traffic light. Hence, in the preparation phase, parametric NLP solvers can be created offline for each possible number of vehicles in a movement, up to a specified maximum. During operation, the parameters are measured online and then utilized with the corresponding pre-generated solver to obtain a new optimal solution for the system. Such an offline preparation phase can be easily accommodated using available software packages such as CasADi \citep{casadi} or ACADO \citep{acado}. The resulting solvers can be reused with updated vehicle states and traffic light timings, so that constructing the optimization problems is carried out in the preparation phase and the online computation is devoted to solving them for the current parameters. The appropriate prepared solver is selected for the current number of vehicles, and the solution and active constraints are determined online. Therefore, changes in the number of vehicles in each movement, and hence in the optimization problem and the corresponding active set, are accommodated by the proposed solution method.
\subsubsection{Distributed optimization}

Using established distributed optimization algorithms \citep{boyd11}, each of the OCPs in the form of \eqref{eq:min_decompose3} can be decomposed into smaller OCPs suitable for parallel in-vehicle implementation. While the development of a distributed control algorithm to solve \eqref{eq:min_decompose3} is not within the scope of this article, efficient algorithms from previous studies for distributed control of vehicles at traffic intersections can be adapted for this purpose \citep{Shi2018,katriniok2021}.

\subsubsection{Total computation time}\label{sec:total-computation}
The total computation time for solving the complete two-level optimization problem at each control update accounts for solving the OCPs in \eqref{eq:min_decompose3}, calculating the derivatives of $\Phi_j(T_j)$, and solving the QP in \eqref{eq:QP} over all SQP iterations, including the evaluations required by line search. The repeated OCP solutions can benefit from the parallel computation and offline preparation described above, while the derivative calculations can reuse available KKT factorizations as explained in Section \ref{sec:solution}.

The decision variables of the upper-level problem \eqref{eq:min_decompose1} consist only of the traffic light timings. Each of the eight movements has a green start time and a green end time in each of the $L$ traffic light cycles. Hence, $\bm{T}$ contains $16L$ entries for the eight movements, independently of the number of vehicles. For a fixed number of cycles, higher traffic demand therefore increases the number and size of the OCPs in \eqref{eq:min_decompose3} while preserving the dimension of the QP in \eqref{eq:QP}. Moreover, the constraints of \eqref{eq:min_decompose1} are linear, and the Hessian of its Lagrange function is block diagonal, with blocks given by \eqref{eq:hessian}. The decomposition thus retains a traffic light timing allocation problem of fixed dimension as the number of vehicles increases, while allowing the vehicle trajectory optimization problems to be solved in parallel.

\begin{remark}
Although experimental deployment of the algorithm is beyond the scope of this paper, sensitivity calculations and techniques similar to those suggested for alleviating computational complexity have been tested and validated in experimental studies \citep{Hult2020Experimental}. In particular, offline preparation, real-time iteration (RTI)-based methods \citep{Diehl2005RTI}, and distributed optimization \citep{Nedic2018Distributed} are conventional steps towards the real-time realizability of optimization-based algorithms and have demonstrated their effectiveness in reducing the online computational effort. Therefore, for practical deployment, established RTI-based methods or distributed optimization techniques could be adopted to address the real-time computational requirements of the proposed method.
\end{remark}

\begin{figure}[t!]
\centering{
\subfloat[\label{fig:scores-demanda}]{

\includegraphics{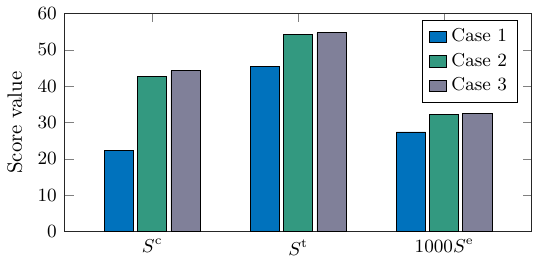}
}\\
\subfloat[\label{fig:scores-demandb}]{
\includegraphics{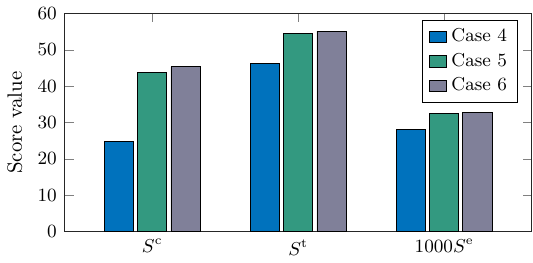}
}\\
\subfloat[\label{fig:scores-demandc}]{
\includegraphics{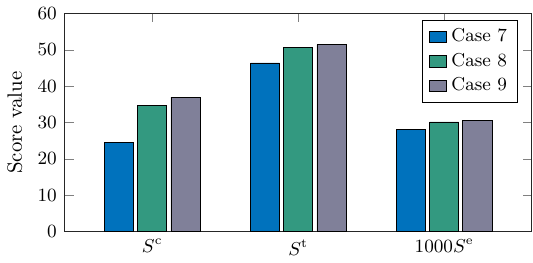}
}
\caption{Comparison of the scores in cases with Demand Profile I (a), Demand Profile II (b), and Demand Profile III (c).}
\label{fig:scores-demand}
}
\end{figure}
\begin{figure}[t!]
\centering{
\subfloat[]{
\includegraphics{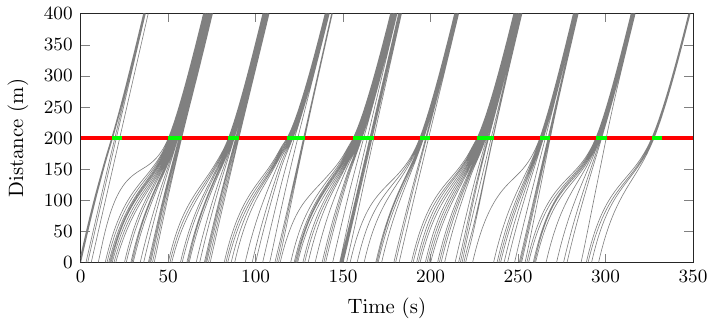}
}\\
\subfloat[]{
\includegraphics{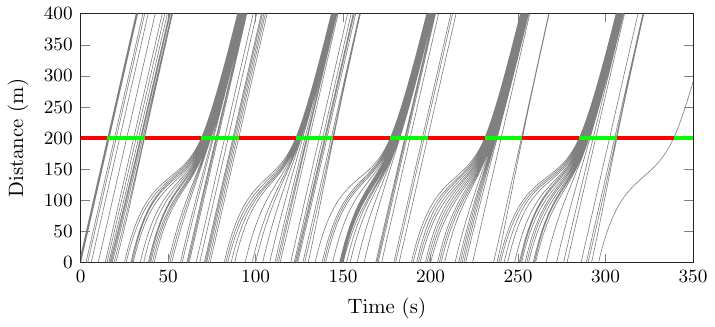}
}\\
\subfloat[]{
\includegraphics{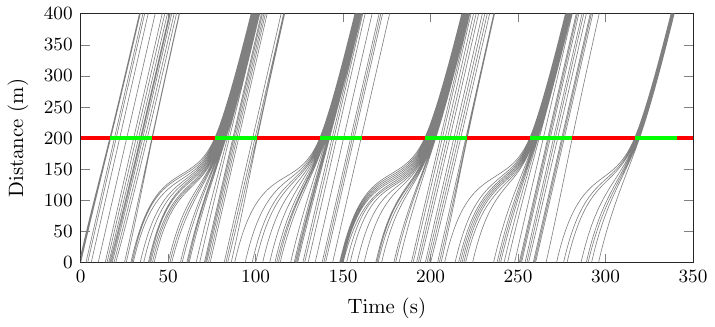}
}
\caption{Vehicle trajectories and traffic light timing in movement 2 for Case 7 (a), Case 8 (b), and Case 9 (c).}
\label{fig:trajectory}
}
\end{figure}
\subsection{Simulation results}

In this section, we present the results obtained through the application of the proposed scheme across various scenarios. To illustrate the viability of this scheme, we consider multiple demand profiles and compare its performance to the case where vehicle trajectory optimization is carried out alongside a traffic light featuring a fixed cycle but demand-dependent timing. Additionally, we provide results that analyze how the scheme's performance is affected by the minimum green time and by the size of the control zone. Finally, we compare the performance of the scheme with the state-of-the-art cooperative traffic light and vehicle speed control method proposed by \citet{Biao_xu-2019}.

\subsubsection{Impact of the traffic demand} 
To assess the proposed cooperative method under various traffic conditions, we conducted simulations with different demand profiles. As a reference for comparison, we benchmarked the proposed control algorithm against a baseline approach in which vehicle trajectories are optimized using \eqref{eq:min_decompose2}, while traffic light timing is not coordinated with the vehicles but is demand-dependent instead. We considered three demand profiles, as detailed in Table \ref{table:demand}, resulting in a total of nine cases. In three of these cases, we adopted the proposed control approach, while the remaining six cases were used for comparison, employing the baseline control approach with various cycle lengths. Table \ref{table:scenarios} summarizes the characteristics of each case, including the control approach, demand profile, and cycle length. 

In each movement, vehicle arrivals are sampled from a Poisson distribution with a rate corresponding to the specific movement's demand. Additionally, the initial speeds of vehicles are randomly selected from a uniform distribution within the range of [$v^\mathrm{d} - 5, v^\mathrm{{max}}$]. If necessary, the arriving speed is adjusted to ensure that a rear-end collision can be avoided when the preceding vehicle brakes with its maximum deceleration. For a fair comparison, cases with the same demand profile share the same initialization parameters. This means that the arrival times and initial speeds of vehicles entering the control zone are identical across these cases.

The scores corresponding to Demand Profiles I, II, and III are compared in Figures \ref{fig:scores-demanda}, \ref{fig:scores-demandb}, and \ref{fig:scores-demandc}, respectively. When evaluating the performance metric $S^\mathrm{c}$, a noticeable improvement is observed when comparing Case 1, Case 4, and Case 7 with their counterparts, indicating the significant benefits of optimizing both trajectories and signal timings jointly. For all demand profiles, this improvement is more pronounced than those in $S^\mathrm{t}$ and $S^\mathrm{e}$, as $S^\mathrm{c}$ reflects the objective that is explicitly minimized in \eqref{eq:min_decompose2}. The improvements in $S^\mathrm{t}$ and $S^\mathrm{e}$, although less prominent, are a by-product of the proposed joint optimization framework.

As an example, in Figure \ref{fig:trajectory}, we illustrate the vehicle trajectories in movement 2 for Case 7, Case 8, and Case 9, all of which use Demand Profile III (refer to Table \ref{table:scenarios}). It is evident from the figure that Case 7, with coordinated traffic light control, leads to smoother vehicle trajectories when passing the intersection. Similar improvements are observed in other movements.

\subsubsection{Impact of the minimum green time} 
\begin{figure}[h]
\centering{
\subfloat[\label{fig:scores-tgvary800}]{
\includegraphics{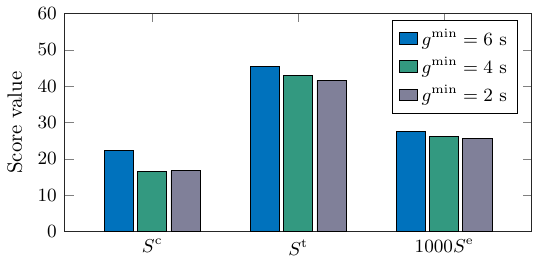}
}\\
\subfloat[\label{fig:scores-tgvary1000}]{
\includegraphics{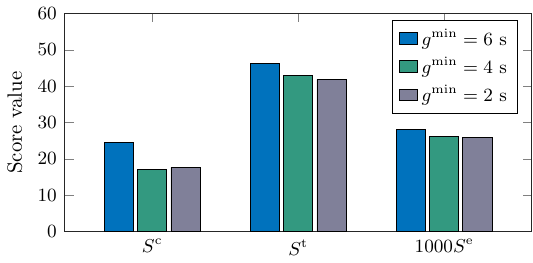}
}
\caption{Impact of the minimum green time on the scores with Demand Profile I (a) and Demand Profile III (b).}
\label{fig:scores-tgvary}

}
\end{figure}

In the proposed algorithm, the minimum green time affects the traffic light timing through \eqref{eq:gmin}. To assess the impact of the minimum green time on the algorithm's performance, we conducted simulations where the minimum green time was chosen from the set $\{2, 4, 6\}$~\si{s}. The results for Demand Profiles I and III are presented in Figures \ref{fig:scores-tgvary800} and \ref{fig:scores-tgvary1000}, respectively. From the figures, it is clear that decreasing the minimum allowed green time leads to performance improvements, although not necessarily in a linear fashion. The most significant improvement is observed when reducing $g^\mathrm{min}$ from 6 to 4~s. Note that, although reducing the minimum green time may generally improve vehicle performance, it is essential to consider other modes of transport, such as pedestrians and cyclists, when selecting a value for $g^\mathrm{min}$, as reducing it may affect these modes differently.

\subsubsection{Impact of the size of the control zone} 
\begin{figure}[h]
\centering{
\includegraphics{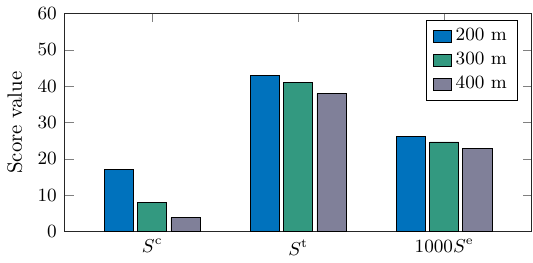}
\caption{Impact of the size of the control zone on the scores. Demand Profile III is used to generate the simulation results.}
\label{fig:communication}
}
\end{figure}

To further evaluate the performance of the proposed scheme under various conditions, we conducted simulations where we adjusted the size of the zone in which vehicle speeds are optimized. In this setup, we varied the radius of the control zone from 200 meters to 400 meters, while keeping all other variables and parameters constant. We selected Demand Profile III for this test. To ensure comparability, we calculated scores over a 200-meter range from the intersection. Figure \ref{fig:communication} illustrates the results, demonstrating that a larger control zone leads to improved performance across all three scores. A longer distance to the intersection allows vehicles to approach more smoothly, which significantly improves the scores. We emphasize that the choice of an appropriate size of the control zone is a tradeoff between performance, communication cost, and equipment availability.

\subsubsection{Comparison with the state-of-the-art} 

We compared the performance of the proposed method with the state-of-the-art (SOTA) method proposed by \citet{Biao_xu-2019}, which utilizes a cooperative approach to traffic light timing optimization and vehicle speed control for automated vehicles in two parts, namely traffic optimization and vehicle speed control. The traffic optimization part uses the enumeration method to optimize the traffic signal timing and vehicle arrival times for a fixed cycle length to minimize the total travel time, whereas the vehicle speed control part uses a pseudospectral method to optimize the vehicle trajectories to minimize vehicle fuel consumption. We implemented the method for the problem formulation given in Section \ref{sec:problem}, and we used the quadratic objective \eqref{eq:cost_fuel} as the objective function for the vehicle speed control part for a fair comparison of the two methods. We refer readers to \citet{Biao_xu-2019} for more details about the SOTA method.

\begin{table*}
\centering
\caption{Scores of the proposed method and the state-of-the-art (SOTA) method for control zones of 200, 300, and 400 meters with Demand Profile III.} 
\label{table:comparison}
\setlength{\tabcolsep}{4pt} 
\begin{tabularx}{\textwidth}{@{} l *{9}{X} @{}} 
\toprule
\multirow{2}{*}{Method} & \multicolumn{3}{c}{200-meter Zone} & \multicolumn{3}{c}{300-meter Zone} & \multicolumn{3}{c}{400-meter Zone} \\
\cmidrule(lr){2-4} \cmidrule(lr){5-7} \cmidrule(lr){8-10}
 & $S^\mathrm{c}$ & $S^\mathrm{t}$ & $1000S^\mathrm{e}$ & $S^\mathrm{c}$ & $S^\mathrm{t}$ & $1000S^\mathrm{e}$ & $S^\mathrm{c}$ & $S^\mathrm{t}$ & $1000S^\mathrm{e}$ \\
\midrule
Proposed & 14.93 & 23.41 & 13.62 & 9.67 & 31.17 & 18.10 & 5.94 & 37.99 & 22.39 \\
SOTA     & 60.12 & 42.93 & 22.64 & 47.32 & 49.26 & 26.75 & 43.78 & 58.71 & 32.40 \\
\bottomrule
\end{tabularx}
\end{table*}

We tested the performance of the two methods using Demand Profile III for control zones of 200, 300, and 400 meters. We used numerical experiments to tune the parameters of the SOTA method to achieve the best performance in each scenario. Based on our findings, we set the fixed cycle length of the SOTA method to 44, 48, and 50 seconds for the control zones of 200, 300, and 400 meters, respectively. Since the SOTA method does not calculate vehicle trajectories beyond the stop line, we utilized the trajectories of the vehicles from the moment they entered the control zone until they passed the stop line for calculating the scores. The corresponding scores of the methods are presented in Table \ref{table:comparison}. The results show that the proposed method outperforms the SOTA method, achieving better performance, total travel time, and explicit fuel consumption scores in each scenario. The superiority of our method over the SOTA method can be explained by two main reasons:

First, the SOTA method calculates the optimal traffic light timings and the arrival times at the end of each cycle with a fixed cycle length, whereas the proposed method calculates the optimal traffic light timings every $t_\mathrm{s}$ seconds with a flexible cycle length. Thus, the proposed method can adjust the traffic light timings at a higher frequency and update the planned trajectories of the vehicles to minimize arrival times and move vehicles at the desired speed, resulting in smaller total travel time and performance scores than the SOTA method.

Second, in our comparison, the SOTA method minimizes acceleration only in the vehicle speed control part, whereas the traffic optimization part calculates optimal traffic light timings to minimize the arrival times without considering acceleration. In contrast, the proposed method integrates traffic light timing optimization and vehicle control, adjusting the traffic light timing at the upper level to reduce acceleration at the lower level. Therefore, it can achieve lower explicit fuel consumption and performance scores than the SOTA method.

\section{Conclusion}\label{sec:conclusion}
In this study, we introduced a two-level cooperative vehicle and traffic light control approach for a single signalized intersection. Our controller is designed to be versatile, not relying on specific vehicle dynamics or objective functions. This flexibility allows us to incorporate complex vehicle models and various objective functions into the framework. Furthermore, we do not constrain traffic light cycle times or require green light durations to be integer values. We also make no assumptions about vehicle arrival speeds and do not simplify vehicle trajectories.

We demonstrated that our optimization problem can be effectively solved using the proposed decomposition method and parametric optimization tools. The algorithm yielded promising results with respect to both the selected objective function and general metrics such as travel time and fuel consumption, especially under high vehicle arrival rates.

Future work will focus on extending the proposed framework to coordinate vehicle trajectories and traffic lights across a network of multiple intersections, with the aim of further improving network efficiency and reducing fuel consumption. Such an extension will need to account for the interactions between neighboring intersections, where vehicle departures influence downstream arrivals and downstream queues may restrict upstream vehicle movement. These interactions would couple the traffic light timing and vehicle trajectory optimization problems across intersections, increasing the computational demand of solving them using a centralized approach. A potential approach is to retain the two-level decomposition at each intersection and to add another level to the hierarchy that acts as an advisor coordinating the traffic lights at neighboring intersections. Exchanging predicted vehicle departures and information on downstream queues with this additional level could allow local problems to be solved in parallel for given coordination decisions, while coordinating traffic light timings and vehicle trajectories across intersections.

Future work will also incorporate vulnerable road users, such as cyclists, into the optimization through prediction models of their movements, allowing traffic signal timings and vehicle trajectories to be jointly optimized with the aim of reducing travel times for all road users.

\section*{Declaration of Generative AI Use}
The authors used OpenAI Codex (GPT-6 Astra) and Anthropic Claude (Opus 5.5) for language improvement and formatting.

\section*{Disclosure statement}
The authors report there are no competing interests to declare.

\section*{Funding}
This work has received funding from the European Union's Horizon 2020 research and innovation programme under the Marie Sk\l{}odowska-Curie grant agreement No 845948. This work has also received funding from the European Research Council (ERC) under the European Union's Horizon 2020 research and innovation programme (Grant agreement No. 101018826 -- ERC Advanced Grant CLariNet).

\section*{Data availability statement}
The data and code that support the findings of this study are available from the corresponding author upon reasonable request.

\bibliographystyle{apacite}
\bibliography{references}

\appendix
\section{}\label{app:sumlog}
To derive \eqref{eq:exp_first_derivative} and \eqref{eq:exp_second_derivative}, let us first define $E$ as
\begin{equation}\label{eq:app_exp}
E = \sum_{b=0}^{N_j} \exp(-\beta \tilde{\Phi}_{j_b}(T_j)).
\end{equation}
The first and second derivatives of $E$ with respect to ${T^o_j}$ are given by

\begin{align}
\dv{E}{T^o_j} &= -\beta\sum_{b=0}^{N_j} \dv{\tilde{\Phi}_{j_b}}{T^o_j}\exp\left(-\beta \tilde{\Phi}_{j_b}(T_j)\right)\\
\dv[2]{E}{{T^o_j}} &= -\beta\sum_{b=0}^{N_j} \left(\dv[2]{\tilde{\Phi}_{j_b}}{{T^o_j}}-\beta \left(\dv{\tilde{\Phi}_{j_b}}{T^o_j}\right)^2\right)\exp(-\beta \tilde{\Phi}_{j_b}(T_j))\label{eq:app_dexp}
\end{align}
Using \eqref{eq:logapproximation} and \eqref{eq:app_exp}--\eqref{eq:app_dexp}, we can write

\begin{align}
\dv{\Phi(T_j)}{T^o_j} &= \frac{-1}{\beta} \displaystyle \frac{\displaystyle\dv{E}{T^o_j}}{E} \nonumber\\
&= \displaystyle\frac{\displaystyle \sum_{b=0}^{N_j} \displaystyle   \dv{\tilde{\Phi}_{j_b}(T_{j})}{T^o_j}  \exp\left(-\beta\tilde{\Phi}_{j_b}(T_{j})\right)}{{\displaystyle\sum_{b=0}^{N_j}  \exp\left(-\beta\tilde{\Phi}_{j_b}(T_{j})\right)}}
\end{align}
and
\begin{align}
&\dv[2]{\Phi(T_j)}{{T^o_j}} = \frac{-1}{\beta} \left( \displaystyle\frac{E\displaystyle\dv[2]{E}{{T^o_j}} - \left( \displaystyle\dv{E}{{T^o_j}}\right)^2}{E^2} \right)\nonumber\\
& =  \frac{1}{\beta}\left( \frac{\displaystyle\dv{E}{T^o_j}}{E} \right)^2 - \frac{1}{\beta} \left( \displaystyle\frac{\displaystyle\dv[2]{E}{{T^o_j}}}{E} \right) \\
& = \displaystyle \beta \left( \frac{\displaystyle \sum_{b=0}^{N_j} \displaystyle \dv{\tilde{\Phi}_{j_b}(T_{j})}{T^o_j} \exp \left(-\beta\tilde{\Phi}_{j_b}(T_{j})\right)} {\displaystyle\sum_{b=0}^{N_j}  \exp\left(-\beta\tilde{\Phi}_{j_b}(T_{j})\right)}   \right)^2  \nonumber\\
 & +  \frac{ \displaystyle \sum_{b=0}^{N_j}\! \exp \left(-\beta\tilde{\Phi}_{j_b}(T_{j})\right) \!\left(\! \dv[2]{\tilde{\Phi}_{j_b}(T_{j})}{{T^o_j}} -\beta \left(\dv{\tilde{\Phi}_{j_b}(T_{j})}{T^o_j}\right)^2\!\right)}             {\displaystyle\sum_{b=0}^{N_j}  \exp \left(-\beta\tilde{\Phi}_{j_b}(T_{j})\right)}.
\end{align}
\section{}\label{app:KKT}
Consider an optimization problem of the general form 
\begin{subequations}\label{eq_app_kkt}
\begin{align}
\min_{u} \quad & \tilde{\Phi}(\tilde{u})\\
\text{s.t.}\quad& \tilde{g}(\tilde{u}) = 0\\
&\tilde{h}(\tilde{u})\leq 0.
\end{align}
\end{subequations}
Let us define $\tilde{z}= [\tilde{u}^{\top},  \tilde{\lambda}^{\top}, \tilde{\mu}_{\mathbb{A}}^{\top}]$ where $\tilde{\lambda}$ is the vector of multipliers associated with the equality constraints $\tilde{g}(\tilde{u})$ and $\tilde{\mu}_{\mathbb{A}}$ denotes the multipliers associated with the strictly active constraints $\tilde{h}_{\mathbb{A}}(\tilde{u})$, i.e., at the solution,  $\tilde{h}_{\mathbb{A}}^i(\tilde{u}) = 0$ and  $\tilde{\mu}^i_{\mathbb{A}} >0$. The KKT optimality conditions state that if the LICQ \citep{NumericalOpt} holds, the primal-dual solution $\tilde{z}^*$ will satisfy the following equation:
\begin{equation}K(\tilde{z}^*) = 
\begin{bmatrix}
\nabla_u\tilde{\mathcal{L}}(\tilde{z}^*)\\
\tilde{g}(\tilde{u}^*)\\
\tilde{h}_{\mathbb{A}}(\tilde{u}^*)
\end{bmatrix}
=0,
\end{equation}
where $\tilde{\mathcal{L}}(\tilde{z})$ is the Lagrange function of problem \eqref{eq_app_kkt}, which is given by
\begin{equation}
\tilde{\mathcal{L}}(\tilde{z}) = \tilde{\Phi}(\tilde{u}) + \tilde{\lambda}^{\top} \tilde{g}(\tilde{u}) + \tilde{\mu}_{\mathbb{A}}^{\top}\tilde{h}_{\mathbb{A}}(\tilde{u}).
\end{equation}

\end{document}